\documentclass[12pt,a4paper]{amsart}
\usepackage{amssymb,latexsym,graphicx}
\usepackage{color}
\usepackage[cp1252]{inputenc} 
\usepackage[T1]{fontenc}
\usepackage{lmodern}
\usepackage[english]{babel}
\usepackage[section]{placeins}
\usepackage[pagebackref=true]{hyperref}
\hypersetup{pdftex,colorlinks=true,linkcolor=blue,citecolor=blue,urlcolor=blue}
\usepackage{version}
\usepackage{currfile}

\excludeversion{old}

\newtheorem{theorem}{Theorem}[section]
\newtheorem{lemma}[theorem]{Lemma}
\newtheorem{proposition}[theorem]{Proposition}

\newtheorem{remark}[theorem]{Remark}

\newtheorem{state}[theorem]{Statement}

\numberwithin{equation}{section}
\numberwithin{figure}{section}
\numberwithin{table}{section}

\definecolor{purple}{RGB}{127,0,255}

\newcommand{\E}{\mathbb{E}}
\newcommand{\N}{\mathbb{N}}

\newcommand{\R}{\mathbb{R}}
\newcommand{\bS}{\mathbb{S}}
\newcommand{\T}{\mathbb{T}}
\newcommand{\Z}{\mathbb{Z}}

\newcommand{\cC}{\mathcal{C}}
\newcommand{\cE}{\mathcal{E}}
\newcommand{\cD}{\mathcal{D}}
\newcommand{\cH}{\mathcal{H}}
\newcommand{\cL}{\mathcal{L}}
\newcommand{\cP}{\mathcal{P}}
\newcommand{\cR}{\mathcal{R}}

\newcommand{\cT}{\mathcal{T}}
\newcommand{\cZ}{\mathcal{Z}}

\newcommand{\mf}{\mathfrak}

\newcommand{\hlambda}{\hat{\lambda}}

\newcommand{\ps}[1]{\langle #1 \rangle}
\newcommand{\alphav}{\alpha{}^{\vee}}

\newcommand{\ecp}{\mathrm{\textsc{ECP}}}

\begin{document}

\title[Courant nodal domain property]{On Courant's nodal domain property for linear combinations of eigenfunctions, Part~{I}}

\author[P. B\'{e}rard]{Pierre B\'erard}
\author[B. Helffer]{Bernard Helffer}

\address{PB: Institut Fourier, Universit\'{e} Grenoble Alpes and CNRS, B.P.74\\ F38402 Saint Martin d'H\`{e}res Cedex, France.}
\email{pierrehberard@gmail.com}

\address{BH: Laboratoire Jean Leray, Universit\'{e} de Nantes and CNRS\\
F44322 Nantes Cedex, France, and LMO, Universit\'e Paris-Sud.}
\email{Bernard.Helffer@univ-nantes.fr}

\date{Revised Oct. 08, 2018  ~(\currfilename)}



\keywords{Eigenfunction, Nodal domain, Courant nodal domain theorem.}

\subjclass[2010]{35P99, 35Q99, 58J50.}

\begin{abstract}
According to Courant's theorem, an eigenfunction as\-sociated with the $n$-th eigenvalue $\lambda_n$ has at most $n$ nodal domains. A footnote in the book of Courant and Hilbert, states that the same assertion is true for any linear combination of eigenfunctions associated with eigenvalues less than or equal to $\lambda_n$. We call this assertion the \emph{Extended Courant Property}.\smallskip

In this paper, we propose simple and explicit examples for which the extended Courant property is false: convex domains in $\R^n$ (hypercube and equilateral triangle), domains with cracks in $\mathbb{R}^2$, on the round sphere $\mathbb{S}^2$, and on a flat torus $\mathbb{T}^2$.
\end{abstract}%

\maketitle

\vspace{1cm}
\begin{center}
To appear in Documenta Mathematica.
\end{center}%
\vspace{1cm}

\section[Introduction]{Introduction}\label{S-intro}
Let $\Omega \subset \R^d$ be a bounded open domain or, more generally, a compact Riemannian manifold with boundary. \medskip

Consider the eigenvalue problem
\begin{equation}\label{E-intro-2}
\left\{
\begin{array}{lll}
- \Delta u &= \lambda\, u &\text{in~} \Omega\,,\\[5pt]
\mf{b}(u) &= 0 &\text{on~} \partial \Omega\,,
\end{array}
\right.
\end{equation}
where $\mf{b}(u)$ is some homogeneous boundary condition on $\partial \Omega$, so that we have a self-adjoint boundary value problem (including the empty condition if $\Omega$ is a closed manifold). For example, we can choose $\mf{d}(u)=u|_{\partial \Omega}$ for the Dirichlet boundary condition, or $\mf{n}(u)=\frac{\partial u}{\partial \nu}|_{\partial \Omega}$ for the Neumann boundary condition.\medskip

Call $H(\Omega,\mf{b})$ the associated self-adjoint extension of $-\Delta$, and list its eigenvalues in nondecreasing order, counting multiplicities, and starting with the index $1$, as
\begin{equation}\label{E-intro-4}
0 \le \lambda_1(\Omega,\mf{b}) < \lambda_2(\Omega,\mf{b}) \le \lambda_3(\Omega,\mf{b}) \le \cdots \,,
\end{equation}
with an associated orthonormal basis of eigenfunctions $\{u_j, j\ge 1\}$.

For any eigenvalue $\lambda$ of $(\Omega,\mf{b})$, define the index
\begin{equation}\label{E-intro-6}
\kappa(\Omega,\mf{b},\lambda) = \min \{k ~|~ \lambda_k(\Omega,\mf{b}) = \lambda\}.
\end{equation}

\textbf{Notation.~}  If $\lambda$ is an eigenvalue of $(\Omega,\mf{b})$, we denote by $\cE(\Omega,\mf{b},\lambda)$ the eigenspace associated with the eigenvalue $\lambda$. \smallskip

We skip $\Omega$ or $\mf{b}$ from the notations, whenever the context is clear.\smallskip

Given a real continuous  function $v$ on $\Omega$, define its \emph{nodal set}
\begin{equation}\label{E-intro-8}
\cZ(v) = \overline{ \left\{ x \in \Omega ~|~ v(x)=0 \right\} }\,,
\end{equation}
and call $\beta_0(v)$ the number of connected components of $\Omega\setminus\cZ(v)$ i.e., the number of \emph{nodal domains} of $v$.

\begin{theorem}\label{T-courant}[Courant, 1923]\\
For any nonzero eigenfunction $u$ associated with $\lambda_n(\Omega,\mf{b})$,
\begin{equation}\label{E-intro-10}
\beta_0(u) \le \kappa\big( \lambda_n(\Omega,\mf{b}) \big) \le n\,.
\end{equation}
\end{theorem}%

Courant's nodal domain theorem can be found in \cite[Chap.~V.6]{CH1953}. \medskip

A footnote in \cite[p.~454]{CH1953}  (second footnote in the German original \cite[p.~394]{CH1931})  indicates: \emph{Any linear combination of the first $n$ eigenfunctions divides the domain, by means of its nodes, into no more than $n$ subdo\-mains. See the G\"{o}ttingen dissertation of H.~Herrmann, Beitr\"{a}ge zur Theorie der Eigenwerte und Eigenfunktionen, 1932.}\medskip

For later reference, we write a precise statement. Given $\lambda \ge 0$, denote by $\cL(\Omega,\mf{b},\lambda)$ the space of linear combinations of eigenfunctions of $H(\Omega,\mf{b})$ associated with eigenvalues less than or equal to $\lambda$,
\begin{equation}\label{E-intro-12}
\cL(\Omega,\mf{b},\lambda) = \left\{ \sum_{\lambda_j(\Omega,\mf{b}) \le \lambda} c_j \, u_j ~|~ c_j \in \R, u_j \in  \cE(\Omega,\mf{b},\lambda_j)\right\}\,.
\end{equation}

\begin{state}\label{P-ECP}[Extended Courant Property]\\
Let $v \in \cL\left( \lambda_n(\Omega,\mf{b}) \right)$ be any linear combination of eigenfunctions associ\-a\-ted with the $n$ first eigenvalues of the eigenvalue problem \eqref{E-intro-2}. Then,
\begin{equation}\label{E-ECP}
\beta_0(v) \le \kappa\big( \lambda_n(\Omega,\mf{b}) \big) \le n\,.
\end{equation}
\end{state}%

We call both Statement~\ref{P-ECP}, and Inequality \eqref{E-ECP}, the \emph{Extended Courant Property}, and refer to them as the $\ecp(\Omega)$, or as the $\ecp(\Omega,\mf{b})$ to insist on the boundary condition $\mf{b}$.\medskip

\subsection{Known results and conjectures}\label{SS-intro-1} We begin by recalling previ\-ous\-ly known results, and conjectures.\smallskip

\textbf{1.~} Statement~\ref{P-ECP} is true for a finite interval, with either the Dirichlet or the Neumann boundary conditions, as well as for the periodic boundary conditions. In dimension $1$, one can actually replace the operator $\frac{d^2}{dx^2}$ by a general Sturm-Liouville operator $\frac{d}{dx}\left(K\frac{d}{dx}\right) + L$, where $K > 0$ and $L$ are functions, see \cite{BH-Sturm} and \cite{BH-ecp2} for more details.\medskip

\textbf{2.~} In \cite{Arn1973}, see also \cite{Arn2011,Ku}, Arnold points out that Statement~\ref{P-ECP} is particularly meaningful in relation to Hilbert's 16th problem. Indeed, let $p$ be a homogeneous real polynomial in $(N+1)$ variables, of even degree $n$. When restricted to the sphere $\bS^N$, or equivalently to the real projective space $\R P^N$, $p$ can be written as a sum of spherical harmonics of even degrees less than or equal to $n$, i.e., as a sum of eigenfunctions of the Laplace-Beltrami operator on $\R P^N$ equipped with the round metric $g_0$. Arnold observes that should $\ecp(\R P^N,g_0)$ be true, then the number of connected components of the complement to the algebraic hypersurface $V_n = p^{-1}(0)$ can be bounded from above by
$$
\mathrm{dim}_{\R}H_0(\R P^N\! \setminus \! V_n,\R) \le C_{N+n-2}^N + 1 \hspace{1cm}(1)
$$

The estimate (1) is known to be true\footnote{We are aware of only one reference for a proof, namely J.~Leydold's thesis \cite{Ley2}, partially published in \cite{Ley3}, using real algebraic geometry.} when $N=2$. It is known to be true when $N=3$ and $n=4$, and false when $N=3$ and $n\ge 6$, with counterexamples constructed by O.~Viro \cite{Vir1979}. \medskip

It follows that $\ecp(\R P^N,g_0)$ is true for $N=2$, and false for $N=3$. Arnold also mentions that $\ecp$ is false for a generic metric $g$ on the sphere, but does not provide any precise statement, nor proof.\medskip

\textbf{Remark}. As mentioned above, $\ecp(\R P^3,g_0)$ is true when restricted to linear combinations of spherical harmonics  of degree less than or equal to $4$. Given any $\lambda_0 > 0$, it is easy to construct a surface $(M,h)$ such that $\ecp(M,h)$ is true for linear combinations of eigenfunctions with eigenvalues less than or equal to $\lambda_0$. Indeed, let $(\bS^1,g_a)$ be the circle with length $2a\pi$. The eigenvalues are the numbers $0$ (with multiplicity $1$), and ${n^2}/{a^2}$, for $n\ge 1$ (with multiplicity $2$). Consider the torus $(M,h_a) = (\bS^1,g_1)\times (\bS^1,g_a)$. Fix some $\lambda_0 > 0$. Then, for $a$ small enough, the eigenfunctions of $(M,h_a)$, associated with eigenvalues less than or equal to $\lambda_0$, correspond to eigenfunctions of the first factor $(\bS^1,g_1)$, for which the extended Courant property is true.\medskip

Fix some $\lambda_0 > 0$. Using the Hopf fibration $\bS^3 \to \bS^2$, and letting the length of the fiber tend to zero as in \cite{BeBo1982}, one can find a metric $g_{\lambda_0}$ on $\bS^3$, such that $\ecp(\bS^3,g_{\lambda_0})$ is true when restricted to linear combinations of eigenfunctions associated with eigenvalues less than or equal to $\lambda_0$.\medskip

\textbf{3.~} In \cite{GZ2003}, Gladwell and Zhu investigate Statement~\ref{P-ECP} (which they call the \emph{Courant-Herrmann conjecture}) for domains in $\R^N$, with the Dirichlet boundary condition. For the Euclidean square $\cC$, they show that $\ecp(\cC)$ is true when restricted to linear combinations of eigen\-func\-tions associated with the first 13 eigenvalues. They make the conjec\-tures that $\ecp$ is true for the square, for rectangles and, more generally for convex domains in $\R^2$. \medskip

They also give numerical evidence that the $\ecp$ is false for more complicated domains (rectangles with perturbed boundary). More precisely, they numerically determine the nodal sets of some linear combinations $c_1 u_1 + c_2 u_2$ of the first two Dirichlet eigenfunctions in such domains, and conclude from the numerical computations that there exist domains for which the linear combinations $c_1 u_1 + c_2 u_2$ have up to $5$ nodal domains. Finally, they also make the conjecture that given any integer $m \ge 2$, there exist a domain, and a linear combination $c_1 u_1 + c_2 u_2$, with $m$ nodal domains.\medskip

\textbf{Remark}. Fix any $\lambda_0 > 0$. Then, for $a$ small enough, $\ecp([0,1]\times [0,a])$ is true when restricted to eigenfunctions associated with eigenva\-lues less than or equal to $\lambda_0$. The reasoning is similar to the one used in the preceding remark.\medskip

\textbf{4.~} Finally, we would like to point out that sums of eigenfunc\-tions appear naturally in several contexts. (i)~ Using the Faber-Krahn inequa\-lity, Pleijel improved Courant's estimate for the Dirichlet Laplacian. In the case of the hypercube, an extension to the Neumann Laplacian could be achieved provided the extended Courant property be true, see \cite{Pl}. (ii)~ As far as their vanishing properties are concerned, eigenfunctions behave like polynomials with degree of the order of the square root of the eigenvalue. From this point of view, as pointed out in \cite{JeLe1999}, it is natural to investigate the vanishing properties of sums of eigenfunctions as well. (iii)~ What is the number of nodal domains of a ``typical'' eigenfunction? One can answer this question, in a probabilistic sense, when eigenspaces have large multiplicities. In a more general framework, one can consider sums of eigenfunctions. We refer to \cite{NaSo2016,Riv2018} and the references in these papers. (iv)~ A similar approach can be made in the framework of Hilbert's 16th problem, see the paper \cite{FyLeLu2015} and its bibliography.

\subsection{Main examples and organization of this paper}\label{SS-intro-2} The purpose of the present paper is to provide simple counter\-examples to the \emph{Exten\-ded Courant Property}.\medskip

In this subsection, we briefly describe the main examples given in this paper. Each of them is directly motivated by a result or by a conjecture mentioned in the previous subsection. \smallskip

\textbf{1.~} Let $\cC_n := [0,\pi]^n$ be the hypercube. In Section~\ref{S-hcube}, we show that $\ecp(\cC_n,\mf{d})$ is false for $n \ge 3$, and that $\ecp(\cC_n,\mf{n})$ is false for $n \ge 4$. This provides convex counterexamples to the $\ecp$ in higher dimensions, for both the Dirichlet or the Neumann boundary conditions.\medskip

\textbf{2.~} Let $\cT_e$ denote the equilateral triangle. In Section~\ref{S-TEQ}, we prove that $\ecp(\cT_e)$ is false for both the Dirichlet and the Neumann boundary conditions.  This provides a convex counterexample to the $\ecp$ in dimension $2$, and therefore a counterexample to one of the conjectures in \cite{GZ2003}. The description of the eigenvalues and eigenfunctions of the equilateral triangle is summarized in Appendix~\ref{A-teq}. \medskip

\textbf{Remark}. By perturbing this example, one can show that there exists a family of smooth strictly convex domains $D_a$ in $\R^2$ such that the $\ecp(D_a,\mf{n})$ is false, see \cite{BH-teqa}. We refer to \cite{BH-ecp2} for other counterexamples related to the equilateral triangle. \medskip

\textbf{3.~} In Section~\ref{S-Rect-c}, we use cracks to perturb the rectangle, or the unit disk. We obtain non-convex, yet simply-connected domains of $\R^2$ which are counterexamples to the $\ecp$. Similarly, in Sections~\ref{S-Tor-c} and \ref{S-Sph-c}, we use cracks to perturb a flat $2$-torus, or the round $2$-sphere, and obtain further counterexamples to the $\ecp$.\medskip

In both cases, we can prescribe the number of nodal domains of the linear combination of eigenfunctions under consideration, thus answer\-ing a conjecture in \cite{GZ2003}.\medskip

\textbf{Remark}. By considering domains with cracks, we are able to provide a rigorous proof of the conjecture proposed in Gladwell and Zhu, based on numerical computations for some domains, and to extend it to the case of non-planar surfaces such as $\T^2$ and $\bS^2$. These examples with cracks also contradict other natural conjectures (such as replacing the minimal labeling $\kappa(\lambda)$ in Statement~\ref{P-ECP}, by a maximal labeling), for which the equilateral triangle is not a counterexample anymore.\medskip

\textbf{4.~} Finally, we would like to point out that some of our examples are relevant to the question of counting the number of connected compon\-ents of the complement of a level line of the second Neumann eigenfunc\-tion, see \cite{BaPa2006}.

\subsection*{Acknowledgements} The authors are very much indebted to Virginie Bon\-naillie-No\"{e}l who produced some simulations and pictures at an early stage of their work on this subject. They thank the referee for his comments.

\section{The hypercube}\label{S-hcube}

\subsection{Preparation}\label{SS-hcubep}

Let $\cC_n(\pi) := ]0,\pi[^{\,n}$ be the \emph{hypercube} of dimension $n$, with either the Dirichlet or the Neumann boundary condition on $\partial \cC_n(\pi)$. A point in $\cC_n(\pi)$ is denoted by $x = (x_1,\ldots,x_n)$. \smallskip

A complete set of eigenfunctions of $-\Delta$ for $(\cC_n(\pi),\mathfrak d)$ is given by the functions
\begin{equation}\label{E-hcube-2d}
\prod_{j=1}^n \sin(k_j \, x_j)\,, \text{~with eigenvalue~} \sum_{j=1}^n k_j^2, \text{~~for~} k_j \in \N\!\setminus\!\{0\}\,.
\end{equation}

A complete set of eigenfunctions of $-\Delta$ for $(\cC_n(\pi),\mathfrak n)$ is given by the functions
\begin{equation}\label{E-hcube-2n}
\prod_{j=1}^n \cos(k_j \, x_j)\,, \text{~with eigenvalue~} \sum_{j=1}^n k_j^2, \text{~~for~} k_j \in \N\,.
\end{equation}

\subsection{Hypercube with Dirichlet boundary condition}\label{SS-hcube-d}

In this sec\-tion, we make use of the classical Chebyshev polynomials $U_k(t), k \in \N$, defined by the relation,
\begin{equation*}
\sin\left( (k+1)t\right) = \sin(t) \, U_k\left( \cos(t) \right) \,,
\end{equation*}
and in particular,
\begin{equation*}
U_0(t) = 1, ~~ U_1(t) = 2t, ~~ U_2(t)=4t^2-1\,.
\end{equation*}

The first Dirichlet eigenvalues of $\cC_n(\pi)$ (as points in the spectrum) are listed in the following table, together with their multiplicities, and eigenfunctions.

\begin{table}[htb]
\caption{First Dirichlet eigenvalues of $\cC_n(\pi)$}\label{T-hcube-2d}
   \centering
   \begin{tabular}{|c|c|c|}
   \hline
   Eigenv. & Mult. & Eigenfunctions \\[5pt]
   \hline
   $n$ & $1$ & $\phi_1(x) := \prod_{j=1}^n \sin(x_j)$\\[5pt]
   \hline
   $n+3$ & $n$ & $\phi_1(x) \, U_1\left( \cos(x_i) \right)$,
   for $1 \le i \le n$\\[5pt]
   \hline
  $n+6$ & $\frac{n(n-1)}{2}$ & $\phi_1(x) \, U_1\left( \cos(x_i) \right) \, U_1\left( \cos(x_j) \right)$,  for $1 \le i < j \le n$\\[5pt]
   \hline
  $n+8$ & $n$ & $\phi_1(x) \, U_2\left( \cos(x_i) \right)$, for
  $1 \le i \le n$\\[5pt]
   \hline
   \end{tabular}
\end{table}

For the above eigenvalues, the index defined in \eqref{E-intro-6} is given by,
\begin{equation}\label{E-hcube-4d}
\kappa(n+3)=2, ~~ \kappa(n+6) = n+2, ~~ \kappa(n+8) = \frac{n(n+1)}{2} + 2 \,.
\end{equation}

\begin{figure}
  \centering
  \includegraphics[scale=0.4]{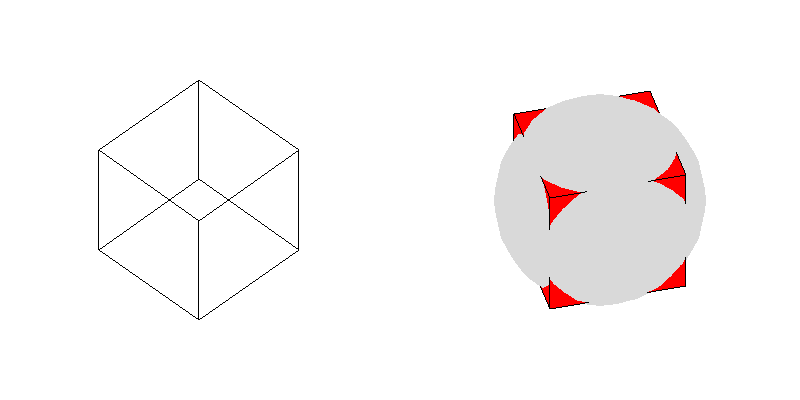}\vspace{-8mm}
  \caption{$3$-dimensional cube}\label{F-hcube-2}
\end{figure}

In order to study the nodal set of the above eigenfunctions or linear combinations thereof, we use the diffeomorphism
\begin{equation}\label{E-hcube-6d}
(x_1,\ldots,x_n) \mapsto \left( \xi_1 = \cos(x_1),\ldots \,, \xi_n = \cos(x_n)  \right) \,,
\end{equation}
from $]0,\pi[^n$ onto $]-1,1[^n$, and factor out the function $\phi_1$ which does not vanish in the open hypercube. We consider the function
$$
\Xi_a(\xi_1,\ldots,\xi_n) = \xi_1^2 + \cdots + \xi_n^2 - a
$$
which corresponds to a linear combination $\Phi_a$ in
$$\cE(\cC_n(\pi),\mf{d},n) \oplus \cE(\cC_n(\pi),\mf{d},n+8).$$
Given some $a$, with $(n-1) < a < n$, the function $\Phi_a$ has $2^n+1$ nodal domains, see Figure~\ref{F-hcube-2} in dimension $3$. For $n \ge 3$, we have $2^n + 1 > \kappa(n+8)$. The function $\Phi_a$ therefore provides a counterexample to the $\ecp$ for the hypercube of dimension at least $3$, with Dirichlet boundary condition.

\begin{proposition}\label{P-hcube-2d}
For $n \ge 3$, the $\ecp(\cC_n(\pi),\mf{d})$ is false.
\end{proposition}%

\textbf{Remark}. An interesting feature of this example is that we get counter\-examples to the $\ecp$ for linear combinations which involve eigenvalues with higher index when $n$ increases. This is also in contrast with the fact that, in dimension $3$, Courant's nodal domain theorem is sharp only for $\delta_1$ and $\delta_2$, \cite{HelKiw2015}.

\subsection{Hypercube with Neumann boundary condition}\label{SS-hcube-n}

In this sec\-tion, we make use of the classical Chebyshev polynomials $T_k(t), k \in \N$, defined by the relation,
\begin{equation*}
\cos(k t) = T_k\left( \cos(t) \right),
\end{equation*}
and in particular,
\begin{equation*}
T_0(t) = 1\,, ~~ T_1(t) = t\,, ~~ T_2(t)=2t^2-1\,.
\end{equation*}

The first Neumann eigenvalues (as points in the spectrum) are listed in the following table, together with their multiplicities, and eigenfunc\-tions.

\begin{table}[htb]
\caption{First Neumann eigenvalues of $\cC_n(\pi)$}\label{T-hcube-2n}
   \centering
   \begin{tabular}{|c|c|c|}
   \hline
   Eigenv. & Mult. & Eigenfunctions \\[5pt]
   \hline
   $0$ & $1$ & $\psi_1(x) := 1$\\[5pt]
   \hline
   $1$ & $n$ & $\cos(x_i)$,
   for $1 \le i \le n$\\[5pt]
   \hline
  $2$ & $\frac{n(n-1)}{2}$ & $\cos(x_i) \, \cos(x_j)$, for $1 \le i < j \le n$\\[5pt]
   \hline
  $3$ & $\frac{n(n-1)(n-2)}{6}$ & $\cos(x_i) \, \cos(x_j)\, \cos(x_k)$, for
  $1 \le i < j < k \le n$\\[5pt]
   \hline
  $4$ & $n+\binom{n}{4}$ & $T_2\left( \cos(x_i) \right)$, for
  $1 \le i \le n$ and \ldots \\[5pt]
   \hline
   \end{tabular}
\end{table}

For these  Neumann eigenvalues, the index defined in \eqref{E-intro-6} is given by,
\begin{equation}\label{E-hcube-4n}
\kappa(2) = n+2\,, \kappa(3) = \frac{n(n+1)}{2} + 2\,, \kappa(4) = \frac{n(n^2+5)}{6}+2\,.
\end{equation}

In order to study the nodal set of the above eigenfunctions or linear combinations thereof, we again use the diffeomorphism \eqref{E-hcube-6d}
and the function $\Xi_a$, which here
corresponds to a linear combination $\Psi_a$ in $\cE(\cC_n(\pi),\mf{n},0) \oplus \cE(\cC_n(\pi),\mf{n},4)$. Given some $a$, with $(n-1) < a < n$, the function $\Psi_a$ has $2^n+1$ nodal domains. For $n \ge 4$, we have $2^n + 1 > \kappa(4)$. The function $\Psi_a$ therefore provides a counterexample to the $\ecp$ for the hypercube of dimension at least $4$, with Neumann boundary condition.

\begin{proposition}\label{P-hcube-2n}
For $n \ge 4\,$, the $\ecp(\cC_n(\pi),\mf{n})$ is false.
\end{proposition}%

\subsection{A stability result for the cube}\label{SS-hcube-3}

According to Subsection~\ref{SS-hcube-d}, the $\ecp(\cC_{3}(\pi),\mf{d})$ is false. Consider the rectangular parallelepiped $\cP_b := ]0,b_1\pi[\times ]0,b_2\pi[\times ]0,b_3\pi[$, with $b=(b_1,b_2,b_3), ~b_i > 0$, and define the $a_i$ by $\sqrt{a_i} \, b_i = 1$. \medskip

The Dirichlet eigenvalues $\delta_i(\cP_b)$ are the numbers $a_1 k_1^2 + a_2 k_2^2 + a_3 k_3^2$, with associated eigenfunctions
\begin{equation}\label{E-hcube-3-4}
\prod_{i=1}^3 \sin \left( k_i \sqrt{a_i} x_i \right), ~~k_i \in \N\!\setminus\!\{0\}.
\end{equation}

The eigenvalues are clearly continuous in the parameters $a_i$. For a generic triple $(a_1,a_2,a_3)$ close enough to $(1,1,1)$, the first $12$ Dirichlet eigenvalues $\delta_i(\cP_b)$ are simple, and correspond to the same type of eigenfunctions as for the ordinary cube (same choices for the triples $(k_1,k_2,k_3)$). This is for example the case if we take $a_1 = 1, a_2 = 1 + \sqrt{2}/100$ and $a_3 = 1 + \sqrt{3}/100$, see the numerical values in Table~\ref{T-hcube-3}, where the Dirichlet eigenvalues are denoted $\delta_i$.\medskip

\begin{table}\label{T-hcube-3}
\caption{Eigenvalues for $(\cC_3(\pi),\mf{d})$ and $(\cP_b,\mf{d})$}
\centering
\begin{tabular}{|c|c|c|c|}
\hline
\text{Index} & \text{Triple} & $\delta_i(\cC_3(\pi))$ & $\delta_i(\cP_b)$\\[5pt]
\hline
$1$ & $(1,1,1)$ & $3$ & $3.016$\\[5pt]
\hline
$2$ & $(2,1,1)$ & $6$ & $6.016$\\[5pt]
\hline
$3$ & $(1,2,1)$ & $6$ & $6.037$\\[5pt]
\hline
$4$ & $(1,1,2)$ & $6$ & $6.042$\\[5pt]
\hline
$5$ & $(2,2,1)$ & $9$ & $9.037$\\[5pt]
\hline
$6$ & $(2,1,2)$ & $9$ & $9.042$\\[5pt]
\hline
$7$ & $(1,2,2)$ & $9$ & $9.063$\\[5pt]
\hline
$8$ & $(3,1,1)$ & $11$ & $11.016$\\[5pt]
\hline
$9$ & $(1,3,1)$ & $11$ & $11.072$\\[5pt]
\hline
$10$ & $(1,1,3)$ & $11$ & $11.085$\\[5pt]
\hline
$11$ & $(2,2,2)$ & $12$ & $12.063$\\[5pt]
\hline
$12$ & $(3,2,1)$ & $14$ & $14.037$\\[5pt]
\hline
\end{tabular}
\end{table}

One can then repeat the arguments of Subsection~\ref{SS-hcube-d}, and conclude that $\ecp(\cC_b,\mf{d})$ is false, so that one has some kind of stability.

\begin{proposition}\label{P-hcube-}
For $b:=(b_1,b_2,b_3)$ close enough to $(1,1,1)$, the $\ecp(\cP_b,\mf{d})$ is false.
\end{proposition}%

Clearly, the same kind of argument can be applied in higher dimension, or for the Neumann boundary condition. \medskip

\textbf{Remark}. Note that the preceding examples still leave open the conjec\-ture made by Gladwell and Zhu that $\ecp$ is true for convex domains in dimension $2$. A counterexample will be given in the next section. \medskip

\section{The equilateral triangle}\label{S-TEQ}

\begin{figure}
  \centering
  \includegraphics[width=7cm]{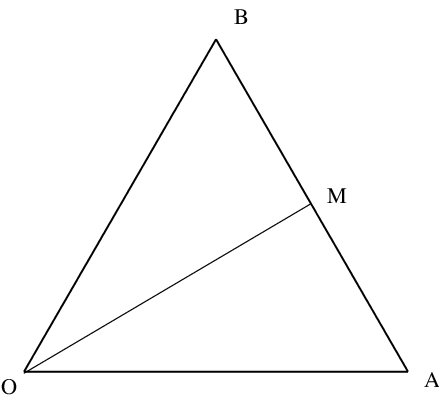}
  \caption{Equilateral triangle $\cT_e = [OAB]$}\label{F-TEQ}
\end{figure}

Let $\cT_e$ denote the equilateral triangle with sides equal to $1$, Figure~\ref{F-TEQ}. The eigenvalues and eigenfunctions of $\cT_e$, with either the Dirichlet or the Neumann condition on the boundary $\partial \cT_e$, can be completely described, see \cite{Be80,McCa2002,McCa2003}, or \cite{BH-lmp16}. We provide a summary in Appendix~\ref{A-teq}. \medskip

In this section, we show that the equilateral triangle provides a counter\-example to the \emph{Extended Courant Property} for both the Dirichlet and the Neumann boundary conditions, contradicting the conjecture of Gladwell and Zhu in dimension $2$.

\subsection{Neumann boundary condition}\label{SS-TEQ-N}

The sequence of Neumann ei\-gen\-values of the equilateral triangle $\cT_e$ begins as follows,
\begin{equation}\label{E-teqN-2}
0 = \lambda_1(\cT_e,\mf{n}) < \frac{16\pi^2}{9} = \lambda_2(\cT_e,\mf{n}) = \lambda_3(\cT_e,\mf{n}) < \lambda_4(\cT_e,\mf{n}) \,.
\end{equation}
The second eigenspace has dimension $2$, and contains one eigenfunction $\varphi^{\mf{n}}_2$ which is invariant under the mirror symmetry with respect to the median $OM$, and another eigenfunction $\varphi^{\mf{n}}_3$ which is anti-invariant under the same mirror symmetry, see Appendix~\ref{A-teq}.\medskip

 More precisely, according to \eqref{E-teq-30}, the function $\varphi^{\mf{n}}_2(x,y)$ can be chosen to be,
\begin{equation}\label{E-teqN-4}
\left\{
\begin{array}{l}
\begin{array}{ll}
\varphi^{\mf{n}}_2(x,y) = &\cos(\frac{4\pi}{3}x) + \cos(\frac{2\pi}{3}(-x+\sqrt{3}y))\\[3pt]
&+ \cos(\frac{2\pi}{3}(x+\sqrt{3}y))\,,
\end{array}%
\end{array}
\right.
\end{equation}
or, more simply,
\begin{equation}\label{E-teqN-6}
\varphi^{\mf{n}}_2(x,y) = 2 \cos\left( \frac{2\pi x}{3} \right) \left( \cos\left( \frac{2\pi x}{3} \right) + \cos\left( \frac{2\pi y}{\sqrt{3}} \right) \right) - 1\,.
\end{equation}

The set $\{\varphi^{\mf{n}}_2 +1 = 0\}$ consists of the two line segments $\{x = \frac{3}{4}\} \cap \cT_e$ and $\{x+\sqrt{3}y=\frac{3}{2}\} \cap \cT_e$, which meet at the point $(\frac{3}{4},\frac{\sqrt{3}}{4})$ on $\partial \cT_e$.\medskip

\begin{figure}
  \centering
  \includegraphics[width=7cm]{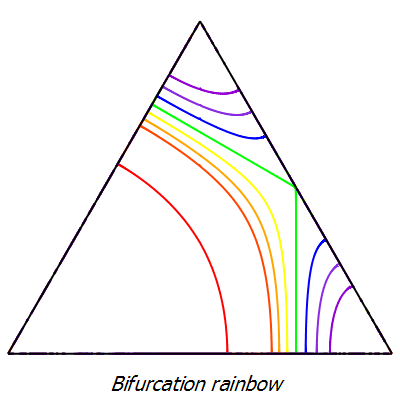}
  \caption{Level sets $\{\varphi^{\mf{n}}_2+a=0\}$ for $ {a\in \{0\,;\, 0.7\,; \, 0.8\,;\,0.9\,;\,1\,;\,1.1\,;\,1.2\,;\,1.3\} }$}\label{F-teqN2}
\end{figure}

The sets $\{\varphi^{\mf{n}}_2 + a=0\}$, with $a \in \{0\,;1-\varepsilon\,;1\,;1+\varepsilon\}$, and small positive $\varepsilon$, are shown in Figure~\ref{F-teqN2}. When $a$ varies from $1-\varepsilon$ to $1+\varepsilon$, the number of nodal domains of $\varphi^{\mf{n}}_2 + a$ in $\cT_e$ jumps from $2$ to $3$, with the jump occurring for $a=1$.\medskip

It follows that $\varphi^{\mf{n}}_2+a=0$, for $1 \le a \le 1.1$, provides a counter\-example to the \emph{Extended Courant Property} for the equilateral triangle with Neumann boundary condition.

\begin{proposition}\label{P-TEQ-2n}
The $\ecp(\cT_e,\mf{n})$ is false.
\end{proposition}%

\textbf{Remark}. The eigenfunction $\varphi^{\mf{n}}_2$ restricted to the hemiequilateral trian\-gle is the second Neumann eigenfunction of $\cT_h = [OAM]$. The restric\-tion of $\varphi^{\mf{n}}_3$ to the hemiequilateral triangle is an eigenfunction of $\cT_h$ with mixed boundary condition (Dirichlet on $OM$ and Neumann on the other sides).

\subsection{Dirichlet boundary condition}\label{SS-TEQ-D}

The sequence of Dirichlet ei\-gen\-values of the equilateral triangle $\cT_e$ begins as follows,

\begin{equation}\label{E-teqD-2}
\lambda_1(\cT_e,\mf{d}) = \frac{16\pi^2}{3} < \lambda_2(\cT_e,\mf{d}) = \lambda_3(\cT_e,\mf{d}) = \frac{112\pi^2}{9} < \lambda_4(\cT_e,\mf{d}).
\end{equation}

 More precisely, according to \eqref{E-teq-36a}, the function $\varphi^{\mf{d}}_1(x,y)$ can be chosen to be,
\begin{equation}\label{E-teqD-4a}
\begin{array}{ll}
\varphi^{\mf{d}}_1(x,y) &=  -8 \, \sin \frac{2\pi y}{\sqrt{3}} \, \sin \pi(x + \frac{y}{\sqrt{3}}) \, \sin \pi(x - \frac{y}{\sqrt{3}})\,,
\end{array}%
\end{equation}
which shows that $\varphi^{\mf{d}}_1$ does not vanish inside $\cT_e\,$.\medskip

The second eigenvalue has multiplicity $2$. It admits one eigenfunction,  $\varphi^{\mf{d}}_2$, which is symmetric with respect to the median $OM$, and given in \eqref{E-teq-40}, and another one, $\varphi^{\mf{d}}_3$, which is anti-symmetric. \medskip

We now consider the linear combination $\varphi^{\mf{d}}_2 + a \, \varphi^{\mf{d}}_1$, with $a$ close to~$1$. The following lemma is the key for reducing the question to the previous analysis.
\begin{lemma}\label{L-teqDN-2}
With the above notation, the following identity holds,
$$
\varphi^{\mf{d}}_2 = \varphi^{\mf{d}}_1 \varphi^{\mf{n}}_2\,.
$$
\end{lemma}

\proof We express the above eigenfunctions in terms of $X :=\cos \frac{2\pi}{3} x$ and $Y := \cos \frac{2\pi}{\sqrt{3}} y$.\medskip

First we observe from \eqref{E-teqN-6} that
$$
\varphi^{\mf{n}}_2(x,y) = 2 X (X+Y) -1\,.
$$
Secondly, we have from \eqref{E-teqD-4a}
$$
\varphi^{\mf{d}}_1(x,y) = 2 \sin \frac{2\pi y}{\sqrt{3}} \,
(8 X^3 - 6 X -2 Y)\,.
$$
Finally, it remains to compute  $\varphi^{\mf{d}}_2$. We start from \eqref{E-teq-40}, and first factorize $\sin \frac{2\pi y}{\sqrt{3}}$ in each line. More precisely, we write,
\begin{equation}\label{E-teqD-6a}
\begin{array}{ll}
 \sin \frac{2\pi}{3} (5 x + \sqrt{3}y) -
 \sin \frac{2\pi}{3} (5 x - \sqrt{3}y) = 2 \sin (\frac{2\pi y}{\sqrt{3}}) \cos ( 5 \frac{2\pi x}{3})\,,  \\[5pt]
 \sin \frac{2\pi}{3} (x - 3\sqrt{3}y) -
 \sin \frac{2\pi}{3} (x + 3\sqrt{3}y) =  -2 \sin (3  \frac{2\pi y}{\sqrt{3}})  \cos (  \frac{2\pi x}{3})\,,   \\[5pt]
 \sin \frac{4\pi}{3} (2 x + \sqrt{3}y) -
 \sin \frac{4\pi}{3} (2 x - \sqrt{3}y) =  2 \sin (2  \frac{2\pi y}{\sqrt{3}})  \cos ( 4  \frac{2\pi x}{3}) \,.
\end{array}%
\end{equation}

We now use the classical Chebyshev polynomials $T_n, U_n$, and the rela\-tions $\cos(n\theta) = T_n(\cos \theta)$ and $\sin(n+1)\theta =\sin(\theta) \, U_n(\cos \theta)$. \medskip

This gives,

\begin{equation*}
\begin{array}{ll}
\varphi^{\mf{d}}_2 & = 2 \sin \frac{2\pi y}{\sqrt{3}} \Big( T_5 (X) - X U_2(Y) + T_4(X) U_1(Y)  \Big)\\[8pt]
&=:  2 \sin \frac{2\pi y}{\sqrt{3}} \, Q (X,Y)  \,.
\end{array}%
\end{equation*}

We find that
$$
Q(X,Y) = 16 X^5 - 20 X^3 + 6 X + 2 Y (8 X^4 - 8 X^2 + 1) - 4 X Y^2\,,
$$
and it turns out that the polynomial $Q(X,Y)$ can be factorized as
$$
 Q(X,Y)= \big( 2X (X+Y) -1 \big) \, (8 X^3 - 6 X - 2 Y)\,,
$$
so that $\varphi^{\mf{d}}_2 = \varphi^{\mf{d}}_1 \varphi^{\mf{n}}_2$.\medskip

In the above computation, we have used the relations,
$$
T_4(X)= 8 X^4-8 X^2+1 \,,\,T_5 (X) = 16 X^5 -20 X^3 + 5 X\,,\,
$$
and
$$
U_1(Y) = 2 Y\,,\, U_2(Y) = 4 Y^2-1 \,.
$$

\hfill \qed

Observing that
$$
\varphi^{\mf{d}}_2 + a \, \varphi^{\mf{d}}_1 = \varphi_1^{\mf{d}} (\varphi^{\mf{n}}_2 +a)\,,
$$
we deduce immediately  from the Neumann result that the function $\varphi^{\mf{d}}_2+a\,\varphi^{\mf{d}}_1$, for $1 \le a \le 1.1$, provides a counter\-example to the \emph{Extended Courant Property} for the equilateral triangle with the Dirichlet bounda\-ry condition.

\begin{proposition}\label{P-TEQ-2d}
The $\ecp(\cT_e,\mf{d})$ is false.
\end{proposition}%

\begin{remark}\label{R-TEQ-6}
Lemma~\ref{L-teqDN-2} is quite puzzling. However, other such identi\-ties do exist. Indeed, consider the square $\cC_2(\pi)$. The first eigenfunction has the form
$$
(x,y) \mapsto \alpha_0 \sin x \sin y \,,
$$
with $\alpha_0\neq 0\,$,
and the second eigenfunctions take the form
$$
(x,y) \mapsto \alpha \sin 2x \sin y + \beta \sin 2y \sin x\,,
$$
with $|\alpha| + |\beta|\neq 0\,$. We can then observe that
$$
\alpha \sin 2x \sin y + \beta \sin 2y \sin x = 2 \sin x \sin y \, (\alpha \cos x + \beta \cos y)\,,
$$
and that $\alpha \cos x + \beta \cos y$ is a Neumann eigenfunction of the square. For $\cC_2(\pi)$, more general relations between Dirichlet and Neumann eigen\-func\-tions follow from the identity $2 T_n = U_n-U_{n-2}$ between Chebyshev polynomials.\\ One can also prove the identity $\varphi_2^{\mf{d}}= a\, \varphi_1^{\mf{d}}\,\varphi_2^{\mf{n}}$ between the eigenfunctions of the right isosceles triangle (for some constant $a$ depending on the normalization of eigenfunctions).
\end{remark}%

\section{Rectangle with a crack}\label{S-Rect-c}

Let $\cR$ be the rectangle $]0,4\pi[\times ]0,2\pi[$. For $0 < a \le 1$, let $C_a := ]0,a]\times \{\pi\}$ and $\cR_a := \cR \setminus C_a$. In this section, we only consider the Neumann boundary condition on $C_a$, and either the Dirichlet or the Neumann boundary condition on $\partial \cR$. The setting is  the one described in \cite[Section~8]{DH1993}.\medskip

We call
\begin{equation}\label{E-RC-2}
\left\{
\begin{array}{l}
0 < \delta_1(0) < \delta_2(0) \le \delta_3(0) \le \cdots\\ \text{resp.}\\[5pt]
0 = \nu_1(0) < \nu_2(0) \le \nu_3(0) \le \cdots
\end{array}%
\right.
\end{equation}
the eigenvalues of $-\Delta$ in $\cR$, with the Dirichlet (resp. the Neumann) boundary condition on $\partial \cR$. They are given by the  numbers $\frac{m^2}{16} + \frac{n^2}{4}$, for pairs $(m,n)$ of positive integers for the Dirichlet problem (resp.  for pairs of non-negative integers for the Neumann problem). Corresponding eigenfunctions are products of sines (Dirichlet) or cosines (Neumann). The eigenvalues are arranged in non-decreasing order, counting multipli\-ci\-ties.\medskip

Similarly, call
\begin{equation}\label{E-RC-2a}
\left\{
\begin{array}{l}
0 < \delta_1(a) < \delta_2(a) \le \delta_3(a) \le \cdots\\  \text{resp.}\\[5pt]
0 = \nu_1(a) < \nu_2(a) \le \nu_3(a) \le \cdots
\end{array}%
\right.
\end{equation}
the eigenvalues of $-\Delta$ in $\cR_a$, with the Dirichlet (resp. the Neumann) boundary condition on $\partial \cR$, and the Neumann boundary condition on $C_a$.\medskip

The first three Dirichlet (resp. Neumann) eigenvalues for the rectangle $\cR$ are as follows.\medskip
\begin{equation}\label{E-RC-6d}
\begin{array}{|c|c|c|c|}
\hline
\text{Eigenvalue} & \text{Value} & \text{Pairs} & \text{Dirichlet eigenfunctions}\\[5pt]
\hline
\delta_1(0) & \frac{5}{16} & (1,1) & \phi_1(x,y) = \sin(\frac{x}{4}) \, \sin(\frac{y}{2})\\[5pt]
\hline
\delta_2(0) & \frac{1}{2} & (2,1) & \phi_2(x,y) = \sin(\frac{x}{2}) \, \sin(\frac{y}{2})\\[5pt]
\hline
\delta_3(0) & \frac{13}{16} & (3,1) & \phi_3(x,y) = \sin(\frac{3x}{4}) \, \sin(\frac{y}{2})\\[5pt]
\hline
\end{array}
\end{equation}
\vspace{5mm}
\begin{equation}\label{E-RC-6n}
\begin{array}{|c|c|c|c|}
\hline
\text{Eigenvalue} & \text{Value} & \text{Pairs} & \text{Neumann eigenfunctions}\\[5pt]
\hline
\nu_1(0) & 0 & (0,0) & \psi_1(x,y) = 1 \\[5pt]
\hline
\nu_2(0) & \frac{1}{16} & (1,0) & \psi_2(x,y) = \cos(\frac{x}{4})\\[5pt]
\hline
\nu_3(0) &  & (0,1) & \psi_3(x,y) =  \cos(\frac{y}{2})\\[5pt]
\nu_4(0) & \frac{1}{4} & (2,0) & \psi_4(x,y) = \cos(\frac{x}{2})\\[5pt]
\hline
\end{array}%
\end{equation}

\begin{figure}
  \centering
  \includegraphics[width=10cm]{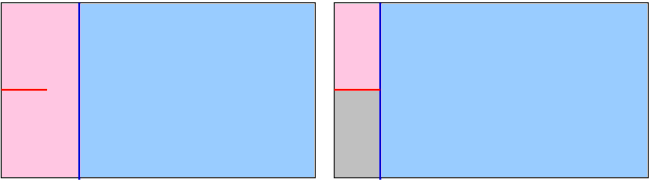}\\
  \caption{Rectangle with a crack (Neumann condition)}\label{F-R1}
\end{figure}

\FloatBarrier

We summarize \cite{DH1993}, Propositions~(8.5), (8.7), (9.5) and (9.9), into the following theorem.

\begin{theorem}[Dauge-Helffer]\label{T-GH}~\\
With the above notation, the following properties hold.
\begin{enumerate}
  \item For $i \ge 1$, the functions $[0,1] \ni a \mapsto \delta_i(a)$, resp. $[0,1] \ni a \mapsto \nu_i(a)$, are non-increasing.
  \item For $i \ge 1$, the functions $]0,1[ \ni a \mapsto \delta_i(a)$, resp. $]0,1[ \ni a \mapsto \nu_i(a)$, are continuous.
  \item For $i \ge 1$, $\lim_{a\to 0+} \delta_i(a) = \delta_i(0)$ and
  $\lim_{a\to 0+} \nu_i(a) = \nu_i(0)$.
\end{enumerate}
\end{theorem}%

It follows that for  $a$ positive, small enough, we have
\begin{equation}\label{E-RC-8}
\left\{
\begin{array}{l}
0 < \delta_1(a) \le \delta_1(0) < \delta_2(a) \le \delta_2(0) < \delta_3(a) \le \delta_3(0)\,, \text{and}\\[5pt]
0 = \nu_1(a) = \nu_1(0) < \nu_2(a) \le \nu_2(0) < \nu_3(a) \le \nu_4(a) \le \nu_3(0)\,.
\end{array}%
\right.
\end{equation}

Observe that for $i = 1$ and $2$,  $\frac{\partial \phi_i}{\partial y}(x,\pi)=0$ and $\frac{\partial \psi_i}{\partial y}(x,y)=0$. It follows that for $a$ small enough, the functions $\phi_1$ and $\phi_2$ (resp. the functions $\psi_1$ and $\psi_2$) are the first two eigenfunctions for $\cR_a$ with the Dirichlet (resp. Neumann) boundary condition on $\partial \cR$, and the Neumann boundary condition on $C_a$, with associated eigenvalues $\frac{5}{16}$ and $\frac{1}{2}$ (resp. $0$ and $\frac{1}{16}$). \medskip

We have
$$
\alpha \phi_1(x,y) + \beta \phi_2(x,y) = \sin(\frac{x}{4})\,\sin(\frac{y}{2})\, \left( \alpha + 2 \beta \cos(\frac{x}{4})\right),
$$
and
$$
\alpha \psi_1(x,y) + \beta \psi_2(x,y) = \alpha + \beta \cos(\frac{x}{4}).
$$
We can choose the coefficients $\alpha, \beta$ in such a way that these linear combinations of the first two eigenfunctions  have two (Figure~\ref{F-R1} left) or three (Figure~\ref{F-R1} right) nodal domains  in $\cR_a\,$. \medskip

\begin{proposition}\label{P-Recc-2}
The $\ecp(\cR_a)$ is false with the Neumann condition on $C_a$, and either the Dirichlet or the Neumann condition on $\partial \cR$.
\end{proposition}%

\begin{remark}\label{R-Recc-2}
In the Neumann case, we can introduce several cracks $\{(x,b_j) ~|~ 0 < x < a_j \}_{j=1}^k$ in such a way that for any $d \in \{2,3, \ldots k+2\}$ there exists a linear combination of $1$ and $\cos(\frac{x}{4})$ with $d$ nodal domains.
\end{remark}%

\begin{remark}\label{R-Recc-4}
Numerical simulations, kindly provided by Virginie Bon\-nail\-lie-No\"{e}l, indicate that the \emph{Extended Courant Property} does not hold for a rectangle with a crack, with the Dirichlet boundary condition on both the boundary of the rectangle, and the crack, \cite{BH-ECP}. Dirichlet cracks appear in another context in \cite{HHOT09} (see also references therein)
\end{remark}%

\begin{figure}[h!bt]
  \centering
  \includegraphics[scale=0.25]{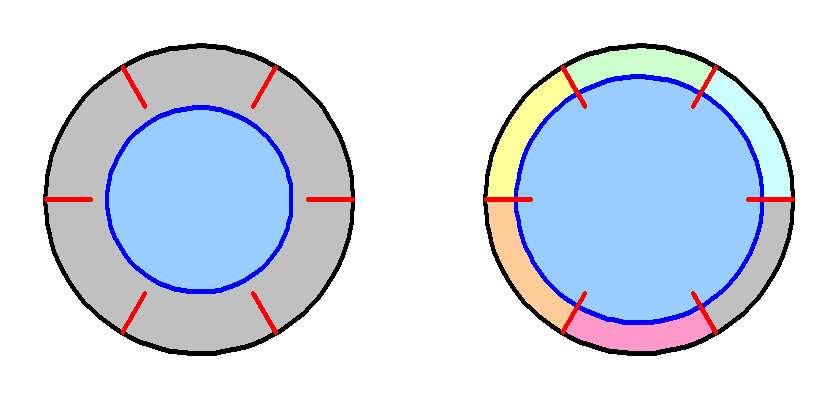}\\
  \caption{Disk with cracks, Neumann condition}\label{F-D1}
\end{figure}

\begin{remark}\label{R-Rect-6}
It is easy to make an analogous construction for the unit disk (Neumann case) with radial cracks. As computed for example in \cite{HPS2016} (Subsection 3.4), the second radial eigenfunction has labelling $6$ ($\lambda_6 \approx 14,68$), and we can introduce six radial cracks to obtain a combination of the two first radial Neumann eigenfunctions with seven nodal domains, see Figure~\ref{F-D1}.
\end{remark}%

\section{The rectangular flat torus with cracks}\label{S-Tor-c}

Consider the flat torus $\T := \R^2/\left( 4\pi \Z \oplus 2\pi Z\right)$. Arrange the eigenvalues in nondecreasing order,
\begin{equation}\label{E-TC-2}
\lambda_1(0) < \lambda_2(0) \le \lambda_3(0) \le \cdots
\end{equation}

The eigenvalues are given by the numbers $\frac{m^2}{4}+n^2$ for $(m,n)$ pairs of integers, with associated complex eigenfunctions
\begin{equation}\label{E-TC-4c}
\exp(im\frac{x}{2})\, \exp(i ny)
\end{equation}
or equivalently,  with real eigenfunctions
\begin{equation}\label{E-TC-4}
\begin{array}{l}
\cos(m\frac{x}{2})\, \cos(ny), \,\cos(m\frac{x}{2})\, \sin(ny),\\[5pt] \sin(m\frac{x}{2})\, \cos(ny), \, \sin(m\frac{x}{2})\, \sin(ny),
\end{array}%
\end{equation}
where $m, n$ are non-negative integers. Accordingly, the first  eigenpairs of $\cT$ are as follows.

\begin{equation}\label{E-RC-6-T}
\begin{array}{|c|c|c|c|}
\hline
\text{Eigenvalue} & \text{Value} & \text{Pairs} & \text{Eigenfunctions}\\[5pt]
\hline
\lambda_1(0) & 0 & (0,0) & \omega_1(x,y) = 1 \\[5pt]
\hline
\lambda_2(0) & & & \omega_2(x,y) = \cos(\frac{x}{2})\\[5pt]
\lambda_3(0) & \frac{1}{4} & (1,0) & \omega_2(x,y) = \sin(\frac{x}{2})\\[5pt]
\hline
\lambda_4(0) &  & & \omega_3(x,y) = \cos(y)\\[5pt]
\lambda_5(0) & 1 & (0,1) & \omega_4(x,y) = \sin(y)\\[5pt]
\hline
\end{array}%
\end{equation}

A typical linear combination of the first three eigenfunctions is of the form $\alpha + \beta \, \sin(\frac{x}{2}-\theta)$\medskip

\begin{figure}
  \centering
  \includegraphics[width=10cm]{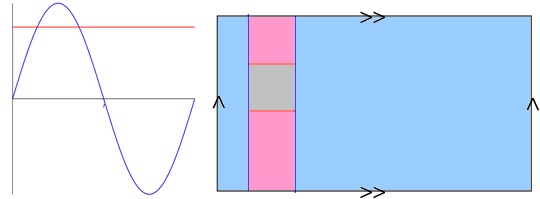}\\
  \caption{Flat torus with two cracks}\label{F-T1}
\end{figure}

Take the torus $\T$, and perform two (or more) cracks parallel to the $x$ axis, and with the same length $a$.  Call $\T_a$ the torus with cracks, see Figure~\ref{F-T1}, and choose the Neumann boundary condition on the cracks. For $a$ small enough, the first three eigenfunctions of the torus $\T$ remain eigenfunctions of the torus with cracks, $\T_a$, with the same $\kappa (\T_a,3)=2$. The  proof is the same as in \cite{DH1993}. We can choose the length $a$ such that the nodal set of $\alpha + \beta \, \sin(\frac{x}{2}-\theta)$ and the two cracks determine three nodal domains.

\begin{proposition}\label{P-Torc-2}
The \emph{Extended Courant Property} is false for the flat torus with cracks (Neumann condition on the cracks).
\end{proposition}%

\section{Sphere $\mathbb S^2$ with cracks}\label{S-Sph-c}

According to \cite[Theorem~1, p.~305]{Ley3}, the $\ecp$ is true on the round sphere $\bS^2$ for sums of spherical harmonics of even (resp. of odd) degree. We consider the geodesic lines  $ z\mapsto (\sqrt{1-z^2} \cos \theta_i, \sqrt{1-z^2}\sin \theta_i, z)$ through the north pole $(0,0,1)$, with distinct $\theta_i \in [0,\pi[$. For example, removing the geodesic segments $\theta_0 = 0$ and $\theta_2= \frac{\pi}{2}$ with $1-z \le a \le 1$, we obtain a sphere $\bS^2_a$ with a crack in the form of a cross. We choose the Neumann boundary condition on the crack.

We can then easily produce a function, in the space generated by the two first eigenspaces of the sphere with a crack, having five nodal domains. \medskip

The first eigenvalue of $\bS^2$ is $\lambda_1(0) = 0$, with corresponding eigenspace of dimension $1$, generated by the function $1$. The next eigenvalues of $\bS^2$ are $\lambda_2(0) = \lambda_3(0) = \lambda_4(0) = 2$ with associated eigenspace of dimension $3$, generated by the functions $x, y, z$. The following eigenvalues of $\bS^2$ are larger than or equal to $6$.

As in \cite{DH1993}, the eigenvalues of $\bS^2_a$ (with Neumann condition on the crack) are non-increasing in $a$, and continuous to the right at $a=0$. More precisely
\begin{equation}\label{E-sph-2}
\left\{
\begin{array}{l}
0 = \lambda_1(a) < \lambda_2(a) \le \lambda_3(a) \le \lambda_4(a) \le 2 < \lambda_5(a) \le 6\,,\\[5pt]
\lim_{a\to 0+}\lambda_i(a) = 2 \text{~for~} i = 2, 3, 4,\\[5pt]
\lim_{a\to 0+}\lambda_5(a) = 6\,.
\end{array}
\right.
\end{equation}

\begin{figure}
  \centering
  \includegraphics[scale=0.5]{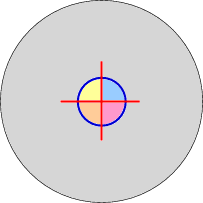}\\
  \caption{Sphere with crack, five nodal domains}\label{F-S1}
\end{figure}

The function $z$ is also an eigenfunction of $\bS^2_a$ with eigenvalue $2$. It follows from \eqref{E-sph-2} that for $a$ small enough, $\lambda_4(a) = 2$, with eigenfunction $z$. For $0 < b < a$, the linear combination $z-b$ has five nodal domains in $\bS^2_a$, see Figure~\ref{F-S1} in spherical coordinates. \medskip

\begin{proposition}\label{P-Sphc-2}
The \emph{Extended Courant Property} is false for the round $2$-sphere with cracks (Neumann condition on the cracks).
\end{proposition}%

\begin{remark}\label{R-Sphc-2}
(1) Removing more geodesic segments around the north pole, we can obtain a linear combination $z-b$ with as many nodal domains as we want. \\
(2) The sphere with cracks, and Dirichlet condition on the cracks, has been considered for another purpose in \cite{HHOT10}.
\end{remark}%

\FloatBarrier

\appendix

\section{Eigenvalues of the equilateral triangle}\label{A-teq}

In this appendix, we recall the description of the eigenvalues of the equilateral triangle. For the reader's convenience, we retain the notation of \cite[Section~2]{BH-lmp16}.\medskip

\subsection{General formulas}

Let $\E^2$ be the Euclidean plane with the canonical orthonormal basis $\{e_1 = (1,0), e_2 = (0,1)\}$, scalar product $\ps{\cdot,\cdot}$ and associated norm $|\cdot|$.\medskip

Consider the vectors
\begin{equation}\label{E-teq-2}
\alpha_1 = (1, - \frac{1}{\sqrt{3}}), \alpha_2 = (0,\frac{2}{\sqrt{3}}), \alpha_3 = (1,\frac{1}{\sqrt{3}}) = \alpha_1 + \alpha_2\,,
\end{equation}
and
\begin{equation}\label{E-teq-4}
\alphav_1 = (\frac{3}{2}, - \frac{\sqrt{3}}{2}), \alphav_2 = (0,\sqrt{3}), \alphav_3 = (\frac{3}{2},\frac{\sqrt{3}}{2}) = \alphav_1 + \alphav_2\,.
\end{equation}

Then
\begin{equation}\label{E-teq-6}
\alphav_i = \frac{3}{2} \alpha_i, |\alpha_i|^2 = \frac{4}{3}, |\alphav_i|^2 = 3.
\end{equation}

Define the mirror symmetries
\begin{equation}\label{E-teq-8}
s_i(x) = x - 2 \frac{\ps{x,\alpha_i}}{\ps{\alpha_i,\alpha_i}}\alpha_i = x - \frac{2}{3}\ps{x,\alphav_i}\alphav_i \,,
\end{equation}
whose axes are the lines
\begin{equation}\label{E-teq-10}
L_i = \{x \in \E^2 ~|~ \ps{x,\alpha_i}=0\}.
\end{equation}

Let $W$ be the group generated by these mirror symmetries. Then,
\begin{equation}\label{E-teq-12}
W = \{1, s_1, s_2, s_3, s_1\circ s_2, s_1\circ s_1\}\,,
\end{equation}
where $s_1 \circ s_2$ (resp. $s_2 \circ s_1$) is the rotation with center the origin and angle $\frac{2\pi}{3}$ (resp. $-\frac{2\pi}{3}$).\medskip

\textbf{Remark}. The above vectors are related to the root system $A_2$ and $W$ is the Weyl group of this root system.\medskip

\begin{figure}[!htb]
	\centering
	\includegraphics[scale=0.35]{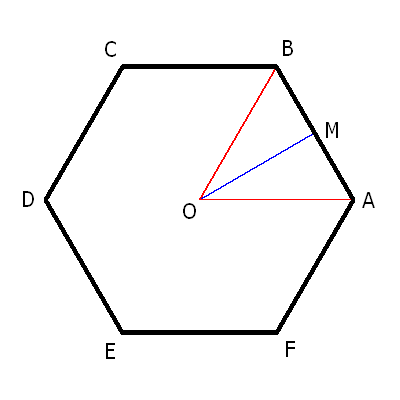}
	\vspace{-8mm}
	\caption{The hexagon $\cH$}\label{F-HT}
\end{figure}

Let
\begin{equation}\label{E-teq-14}
\Gamma = \Z \alphav_1 \oplus \Z \alphav_2
\end{equation}
be the (equilateral) lattice. The set
\begin{equation}\label{E-teq-16}
\cD_{\Gamma} = \{s \alphav_1 + t \alphav_2 ~|~ 0 \le s,t \le 1\}
\end{equation}
is a fundamental domain for the action of $\Gamma$ on $\E^2$. Another fundam\-en\-tal domain is the closure of the open hexagon (see Figure~\ref{F-HT})
\begin{equation}\label{E-teq-18a}
\cH = [A,B,C,D,E,F]\,,
\end{equation}
whose vertices are given by
\begin{equation}\label{E-teq-18b}
\left\{
\begin{array}{l}
A = (1,0) ; B = (\frac{1}{2},\frac{\sqrt{3}}{2}) ; C = (-\frac{1}{2},\frac{\sqrt{3}}{2}) ;\\[5pt]
D = (-1,0) ; E = (-\frac{1}{2},-\frac{\sqrt{3}}{2}) ; F = (\frac{1}{2},-\frac{\sqrt{3}}{2})\,.
\end{array}
\right.
\end{equation}

Call $\cT_e$ the equilateral triangle
\begin{equation}\label{E-teq-18c}
\cT_e = [O,A,B]\,,
\end{equation}
where $O = (0,0)$.

Let $\Gamma^*$ be the dual lattice of the lattice $\Gamma$, defined by
\begin{equation}\label{E-teq-20}
\Gamma^* = \{x \in \E^2 ~|~ \forall \gamma \in \Gamma , \ps{x,\gamma} \in \Z\}\,.
\end{equation}

Then,
\begin{equation}\label{E-teq-20a}
\left\{
\begin{array}{l}
\Gamma^* = \Z \varpi_1 \oplus \Z \varpi_2 \,,\\[5pt]
\text{where~} \varpi_1 = (\frac{2}{3},0) \text{~and~} \varpi_2 = (\frac{1}{3},\frac{1}{\sqrt{3}})\,.
\end{array}
\right.
\end{equation}

Define the set $C$ (an open Weyl chamber of the root system $A_2$),
\begin{equation}\label{E-teq-22}
C = \{x \varpi_1 + y \varpi_2 ~|~ x,y > 0\}\,,
\end{equation}
and let $\T_e$ denote the equilateral torus $\E^2/\Gamma$. \medskip

A complete set of orthogonal (not normalized) eigenfunctions of $-\Delta$ on $\T_e$ is given (in complex form) by the exponentials
\begin{equation}\label{E-teq-24}
\phi_p(x) = \exp (2i\pi \ps{x,p}) \text{~where~} x \in \E^2 \text{~and~} p \in \Gamma^* \,.
\end{equation}
Furthermore, for $p = m \varpi_1 + n \varpi_2$, with $m, n \in \Z$, the multiplicity of the eigenvalue $\hlambda(m,n) = 4\pi^2 |p|^2 = \frac{16\pi^2}{9}(m^2 + mn + n^2)$ is equal to the number of points $(k,\ell)$ in $\Z^2$ such that $k^2+k\ell + \ell^2 = m^2+mn+n^2$.\medskip

The closure of the equilateral triangle $\cT_e$ is a fundamental domain of the action of the semi-direct product $\Gamma \rtimes W$ on $\E^2$ or equivalently, a fundamental domain of the action of $W$ on $\T_e^2$.\medskip

For the following proposition, we refer to \cite{Be80}.

\begin{proposition}\label{P-teq-2}
Complete orthogonal (not normalized) sets of ei\-gen\-functions of the equilateral triangle $\cT_e$ in complex form are given, respectively for the Dirichlet (resp. Neumann) boundary condition on $\partial \cT_e$, as follows.
\begin{enumerate}
  \item Dirichlet boundary condition on $\partial \cT_e$. The family is
  \begin{equation}\label{E-teq-30d}
  \Phi^{\mf{d}}_p(x) = \sum_{w \in W} det(w) exp (2 i \pi \ps{x,w(p)})
  \end{equation}
  with $p \in C \cap \Gamma^*$. Furthermore, for $p = m\varpi_1 + n\varpi_2$, with $m,n$ positive integers, the multiplicity of the eigenvalue $4\pi^2 |p|^2$ is equal to the number of solutions $q \in C \cap \Gamma^*$ of the equation $|q|^2 = |p|^2$.

  \item Neumann boundary condition on $\partial \cT_e$. The family is
  \begin{equation}\label{E-teq-30e}
  \Phi^{\mf{n}}_p(x) = \sum_{w \in W} exp (2 i \pi \ps{x,w(p)})
  \end{equation}
 with $p \in \overline{C} \cap \Gamma^*$. Furthermore, for $p = m\varpi_1 + n\varpi_2$, with $m,n$ non-negative integers, the multiplicity of the eigenvalue $4\pi^2 |p|^2$ is equal to the number of solutions $q \in \overline{C} \cap \Gamma^*$ of the equation $|q|^2 = |p|^2$.
\end{enumerate}
\end{proposition}%

\textbf{Remark}. To obtain corresponding complete orthogonal sets of real eigenfunctions, it suffices to consider the functions
$$
C_p = \Re (\Phi_p) \text{~and~} S_p = \Im (\Phi_p)\,.
$$

For $p = m \varpi_1 + n \varpi_2$, with $m, n \in \N \setminus \{0\}$ for the Dirichlet boundary condition (resp. $m, n \in \N$ for the Neumann boundary condition), we denote these functions by $C_{m,n}$ and $S_{m,n}$.\medskip

In order to give explicit formulas for the first eigenfunctions, we have to examine the action of the group $W$ on the lattice $\Gamma^*$. A simple calculation yields the following table in which we simply denote \break  $m \varpi_1 + n \varpi_2$ by $(m,n)$.

\begin{equation}\label{E-teq-26}
\begin{array}{|c|c|c|c|}
\hline
w &  (m,n) & w(m,n) & \det(w) \\[5pt]
\hline
1 & (m,n) & (m,n) & 1\\[5pt]
\hline
s_1 & (m,n) & (-m,m+n) & -1\\[5pt]
\hline
s_2 & (m,n) & (m+n,-n)& -1 \\[5pt]
\hline
s_3 & (m,n) & (-n,-m)& -1 \\[5pt]
\hline
s_1\circ s_2 & (m,n) & (-m-n,m) & 1\\[5pt]
\hline
s_2\circ s_1 & (m,n) & (n,-m-n) & 1\\[5pt]
\hline
\end{array}%
\end{equation}\medskip

\textbf{Remark}. The above table should be compared with \cite[Table]{BH-lmp16},  in which there is a slight unimportant error (the lines $s_1\circ s_2$ and $s_2 \circ s_1$ are interchanged).\medskip

\textbf{Remark}. Using the above chart, one can easily prove the following relations.\medskip

\begin{equation}\label{E-teq-28}
\left\{
\begin{array}{lll}
C^{\mf{d}}_{n,m} = - C^{\mf{d}}_{m,n} &\text{~and~}& S^{\mf{d}}_{n,m} = S^{\mf{d}}_{m,n}\,,\\[5pt]
C^{\mf{n}}_{n,m} = C^{\mf{n}}_{m,n} &\text{~and~}& S^{\mf{n}}_{n,m} = - S^{\mf{n}}_{m,n}\,.
\end{array}
\right.
\end{equation}\medskip

\subsection{Neumann boundary condition, first three eigenfunctions} The first Neumann eigenvalue of $\cT_e$ is $0$, corresponding to the point $0 = (0,0) \in \Gamma^*$, with first eigenfunction $\varphi_1 \equiv 1$ up to scaling. \medskip

The second Neumann eigenvalue corresponds to the pairs $(1,0)$ and $(0,1)$. According to the preceding remark, it suffices to consider $C^{\mf{n}}_{1,0}$ and $S^{\mf{n}}_{1,0}$. Using Proposition~\ref{P-teq-2}, and the table \eqref{E-teq-26}, we find that, at the point $ [s,t] = s\alphav_1+t\alphav_2$,
\begin{equation}\label{E-teq-29}
\left\{
\begin{array}{l}
C^{\mf{n}}_{1,0}([s,t]) = 2 \left( \cos (2\pi s) + \cos (2\pi (-s+t)) + \cos (2\pi t)\right)\,,\\[5pt]
S^{\mf{n}}_{1,0}([s,t]) = 2 \left( \sin (2\pi s) + \sin (2\pi (-s+t)) -\sin (2\pi t)\right)\,.
\end{array}%
\right.
\end{equation}

Up to a factor $2$, this gives the following two independent eigenfunctions for the Neumann eigenvalue $\frac{16\pi^2}{9}$, in the $(x,y)$ variables, with $(x,y) = \left( \frac{3}{2}s, -\frac{\sqrt{3}}{2}s + \sqrt{3}t \right)$ or $(s,t) = \left( \frac{2}{3}x, \frac{1}{3}x + \frac{1}{\sqrt{3}}y \right)$,
\begin{equation}\label{E-teq-30}
\left\{
\begin{array}{ll}
\varphi^{\mf{n}}_2(x,y) = &\cos(\frac{4\pi}{3}x) + \cos(\frac{2\pi}{3}(-x+\sqrt{3}y))\\[3pt]
&+ \cos(\frac{2\pi}{3}(x+\sqrt{3}y))\,,\\[5pt]
\varphi^{\mf{n}}_3(x,y) = & \sin(\frac{4\pi}{3}x) + \sin(\frac{2\pi}{3}(-x+\sqrt{3}y))\\[3pt]
&- \sin(\frac{2\pi}{3}(x+\sqrt{3}y))\,.
\end{array}%
\right.
\end{equation}

The first eigenfunction is invariant under the mirror symmetry with respect to the median $OM$ of the equilateral triangle, see Figure~\ref{F-TEQ}. The second eigenfunction is anti-invariant under the mirror symmetry with respect to this median. Its nodal set is equal to the median itself.

\subsection{Dirichlet boundary condition, first three eigenfunctions} The first Dirichlet eigenvalue of $\cT_e$ is $\delta_1(\cT_e) = \frac{16\pi^2}{3}$. A first eigenfunction is given by $S^{\mf{d}}_{1,1}$. Using Proposition~\ref{P-teq-2} and Table~\ref{E-teq-26}, we find that this eigenfunction is given, at the point $[s,t]= s\alphav_1+t\alphav_2$, by the formula
\begin{equation}\label{E-teq-36}
\left\{
\begin{array}{ll}
\varphi^{\mf{d}}_1([s,t]) = & 2 \sin 2\pi (s+t) + 2 \sin 2\pi (s-2t) \\[3pt]
& + 2 \sin 2\pi (t-2s)\,.
\end{array}%
\right.
\end{equation}

 Substituting the expressions of $s$ and $t$ in terms of $x$ and $y$, one obtains the formula,
\begin{equation}\label{E-teq-36a}
\begin{array}{ll}
\varphi^{\mf{d}}_1(x,y) = & 2 \sin\left( 2\pi (x + \frac{y}{\sqrt{3}}) \right)
- 2 \sin \left( 4 \pi \frac{y}{\sqrt{3}} \right)\\[3pt]
& - 2 \sin\left( 2\pi (x - \frac{y}{\sqrt{3}}) \right)\,,
\end{array}%
\end{equation}

The second Dirichlet eigenvalue has multiplicity $2$,
$$ \delta_2(\cT_e) = \delta_3(\cT_e) = \frac{112\pi^2}{9}\,.
$$
The eigenfunctions $C^{\mf{d}}_{2,1}$ and $S^{\mf{d}}_{2,1}$ are respectively anti-invariant and invariant under the mirror symmetry with respect to $[OM]$, with values at the point $[(s,t)]$ given by the formulas,

\begin{equation}\label{E-teq-38}
\left\{
\begin{array}{ll}
\varphi^{\mf{d}}_2([s,t]) = & \sin 2\pi(2s+t) + \sin 2\pi(s+2t)\\[3pt]
& +\sin 2\pi(2s-3t) - \sin 2\pi(3s-2t)\\[3pt]
& +\sin 2\pi(s-3t) - \sin 2\pi(3s-t)\,,\\[5pt]
\varphi^{\mf{d}}_3([s,t]) = & \cos 2\pi(2s+t) - \cos 2\pi(s+2t)\\[3pt]
& -\cos 2\pi(2s-3t) + \cos 2\pi(3s-2t)\\[3pt]
& +\cos 2\pi(s-3t) - \cos 2\pi(3s-t)\,.
\end{array}%
\right.
\end{equation}

Substituting the expressions of $s$ and $t$ in terms of $x$ and $y$, one obtains the formulas,

\begin{equation}\label{E-teq-40}
\begin{array}{ll}
\varphi^{\mf{d}}_2(x,y) = & \sin\left( \frac{2\pi}{3} (5 x + \sqrt{3}y) \right)
- \sin\left( \frac{2\pi}{3} (5 x - \sqrt{3}y) \right)\\
& + \sin\left( \frac{2\pi}{3} (x - 3\sqrt{3}y) \right)
- \sin\left( \frac{2\pi}{3} (x + 3\sqrt{3}y) \right)\\
& + \sin\left( \frac{4\pi}{3} (2 x + \sqrt{3}y) \right)
- \sin\left( \frac{4\pi}{3} (2 x - \sqrt{3}y) \right)\,.
\end{array}%
\end{equation}
and
\begin{equation}\label{E-teq-42}
\begin{array}{ll}
\varphi^{\mf{d}}_3(x,y) = & \cos\left( \frac{2\pi}{3} (5 x + \sqrt{3}y) \right)
- \cos\left( \frac{2\pi}{3} (5 x - \sqrt{3}y) \right)\\
& + \cos\left( \frac{2\pi}{3} (x - 3\sqrt{3}y) \right)
- \cos\left( \frac{2\pi}{3} (x + 3\sqrt{3}y) \right)\\
& + \cos\left( \frac{4\pi}{3} (2 x + \sqrt{3}y) \right)
- \cos\left( \frac{4\pi}{3} (2 x - \sqrt{3}y) \right)\,.
\end{array}%
\end{equation}

\bibliographystyle{plain}

\begin{thebibliography}{}

\end{thebibliography}


\begin{thebibliography}{1}


\bibitem{Arn1973} V. Arnold.
\newblock The topology of real algebraic curves (the works of Petrovskii and their development).
\newblock Uspekhi Math. Nauk. 28:5 (1973), 260--262.
\newblock English translation in \cite{Arn2014}.

\bibitem{Arn2011} V. Arnold.
\newblock Topological properties of eigenoscillations in mathematical physics.
\newblock Proc. Steklov Inst. Math. 273 (2011), 25--34.

\bibitem{Arn2014} V. Arnold.
\newblock Topology of real algebraic curves (Works of I.G.~Petrovskii and their development). Translated from \cite{Arn1973} by Oleg Viro.
\newblock In Collected works, Volume II. Hydrodynamics, Bifurcation theory and Algebraic geometry, 1965--1972.
\newblock Edited by A.B.~Givental, B.A.~Khesin, A.N.~Varchenko, V.A.~Vassilev, O.Ya.~Viro. Springer 2014.
\url{http://dx.doi.org/10.1007/978-3-642-31031-7}~.
\newblock Chapter 27, pages 251--254. \url{http://dx.doi.org/10.1007/978-3-642-31031-7_27}~.

\bibitem{BaPa2006} R. Ba{\~{n}}uelos and M. Pang.
\newblock Level sets of Neumann eigenfunctions.
\newblock Indiana University Math. J. 55:3 (2006), 923--939.

\bibitem{Be80} P. B\'{e}rard.
\newblock Spectres et groupes cristallographiques.
\newblock Inventiones Math. 58 (1980), 179--199.

\bibitem{BH-lmp16} P. B\'{e}rard and B. Helffer.
\newblock Courant-sharp eigenvalues for the equilateral torus, and for the equilateral triangle.
\newblock Letters in Math. Physics 106:12 (2016), 1729--1789.

 \bibitem{BH-Sturm} P. B\'{e}rard and B. Helffer.
\newblock Sturm's theorem on zeros of linear combinations of eigenfunctions.
\newblock arXiv:1706.08247.
\newblock To appear in Expositiones Mathematicae.

\bibitem{BH-teqa} P. B\'{e}rard and B. Helffer.
\newblock Level sets of certain Neumann eigenfunctions under deformation of Lipschitz domains. Application to the Extended Courant Property.
\newblock arXiv:1805.01335.

\bibitem{BH-ECP} P. B\'{e}rard and B. Helffer.
\newblock On Courant's nodal domain property for linear combinations of eigenfunctions.
\newblock arXiv:1705.03731v3 (23 Oct 2017).

\bibitem{BH-ecp2} P. B\'{e}rard and B. Helffer.
\newblock On Courant's nodal domain property for linear combinations of eigenfunctions, Part~{II}.
\newblock arXiv:1803.00449 (version $\ge$ 2).

\bibitem{BeBo1982} L.~B\'{e}rard Bergery and J.P.~Bourguignon.
\newblock Laplacians and Riemannian submersions with totally geodesic fibers.
\newblock Ill. J. Math. 26 (1982), 181--200.

\bibitem{CH1931} R. Courant and D. Hilbert.
\newblock Methoden der mathematischen Physik.
\newblock Erster Band. Zweite verbesserte Auflage.
\newblock Julius Springer 1931.

\bibitem{CH1953} R. Courant and D. Hilbert.
\newblock Methods of mathematical physics. Vol. 1.
\newblock First English edition. Interscience, New York 1953.

\bibitem{DH1993} M. Dauge and B. Helffer.
\newblock Eigenvalues variation II. Multidimensional problems.
\newblock J. Diff. Eq. 104 (1993), 263--297.

\bibitem{FyLeLu2015} Y.~Fyodorov, A.~Lerario, and E.~Lundberg.
\newblock On the number of connected components of random algebraic hypersurfaces.
\newblock arXiv:1404.5349.
\newblock J. Geom. Phys. 95 (2015), 1--20.

\bibitem{GZ2003} G. Gladwell and H. Zhu.
\newblock The Courant-Herrmann conjecture.
\newblock ZAMM - Z. Angew. Math. Mech. 83:4 (2003), 275--281.

\bibitem{HHOT09} B. Helffer and T. Hoffmann-Ostenhof and S. Terracini.
\newblock Nodal domains and spectral minimal partitions.
\newblock Ann. Inst. H. Poincar\'{e} Anal. Non Lin\'{e}aire 26 (2009), 101--138.

\bibitem{HHOT10} B. Helffer and T. Hoffmann-Ostenhof and S. Terracini.
\newblock On spectral minimal partitions: the case of the sphere.
\newblock In \emph{Around the Research of Vladimir Maz'ya III}.
\newblock International Math. Series, Springer, Vol. 13, p. 153--178 (2010).

\bibitem{HelKiw2015} B.~Helffer and R.~Kiwan.
\newblock Dirichlet eigenfunctions on the cube, sharpening the Courant nodal inequality.
\newblock arXiv: 1506.05733.
\newblock  In \emph{Functional analysis and operator theory for quantum physics. The Pavel Exner anniversary volume. Dedicated to Pavel Exner on the occasion of his 70th birthday.} Dittrich, Jaroslav (ed.) et al. European Mathematical Society. EMS Series of Congress Reports, 353-371 (2017).

\bibitem{HPS2016} B. Helffer and M. Persson-Sundqvist.
\newblock On nodal domains in Euclidean balls.
\newblock arXiv:1506.04033v2. Proc. Amer. Math. Soc. 144:11 (2016), 4777--4791.

\bibitem{JeLe1999} D. Jerison and G. Lebeau.
\newblock Nodal sets of sums of eigenfunctions.
\newblock In \emph{Harmonic analysis and partial differential equations.
Essays in honor of Alberto Calder\'{o}n}. Edited by M. Christ, C. Kenig and C. Sadosky.
Chicago Lectures in Mathematics, 1999.
\newblock Chap. 14, pp. 223-239.

 \bibitem{Ku} N. Kuznetsov.
\newblock On delusive nodal sets of free oscillations.
\newblock Newsletter of the European Mathematical Society 96 (2015), 34--40.

\bibitem{Ley2} J. Leydold.
\newblock On the number of nodal domains of spherical harmonics.
\newblock PhD Thesis, Vienna University (1992).

\bibitem{Ley3} J. Leydold.
 \newblock On the number of nodal domains of spherical harmonics.
\newblock Topology 35 (1996), 301--321.

\bibitem{McCa2003} J.B. McCartin.
\newblock Eigenstructure of the Equilateral Triangle,
 Part I: The Dirichlet Problem.
\newblock SIAM Review, 45:2 (2003), 267--287.

\bibitem{McCa2002} J.B. McCartin.
\newblock Eigenstructure of the Equilateral Triangle,
 Part II: The Neumann Problem.
\newblock Mathematical Problems in Engineering 8:6 (2002), 517--539.

\bibitem{NaSo2016} F.~Nazarov and M.~Sodin.
\newblock Asymptotic laws for the spacial distribution and the number of connected
components of zero sets of Gaussian random functions.
\newblock Zh. Mat. Fiz. Anal. Geom. 12:3 (2016), 205--278.

\bibitem  {Pl} {\AA}.~Pleijel.
\newblock Remarks on Courant's nodal theorem.
\newblock Comm. Pure. Appl. Math. 9 (1956), 543--550.

\bibitem{Riv2018} A. Rivera.
\newblock Expected number of nodal components for cut-off fractional Gaussian fields.
\newblock arXiv:1801.06999.

\bibitem{Vir1979} O. Viro.
\newblock Construction of multi-component real algebraic surfaces.
\newblock Soviet Math. dokl. 20:5 (1979), 991--995.

\end{thebibliography}

\end{document}